
\parindent=0pt
\def\square{\,\hbox{\vrule\vbox{\hrule\phantom{N}\hrule}\vrule}\,}



%

\def\Cp{\oalign{$C$\cr \noalign{\vskip0.3ex}\hidewidth$\scriptstyle p$\hidewidth}{}}
\def\Cpm{\oalign{$C$\cr \noalign{\vskip0.3ex}\hidewidth$\scriptstyle p-1$\hidewidth}{}}

\def\C0{\oalign{$C$\cr \noalign{\vskip0.3ex}\hidewidth$\scriptstyle 0$\hidewidth}{}}
\def\Ca{\oalign{$C$\cr \noalign{\vskip0.3ex}\hidewidth$\scriptstyle 1$\hidewidth}{}}

\def\Cb{\oalign{$C$\cr \noalign{\vskip0.3ex}\hidewidth$\scriptstyle 2$\hidewidth}{}}

\def\Cc{\oalign{$C$\cr \noalign{\vskip0.3ex}\hidewidth$\scriptstyle 3$\hidewidth}{}}
\def\Cd{\oalign{$C$\cr \noalign{\vskip0.3ex}\hidewidth$\scriptstyle 4$\hidewidth}{}}
\def\Ce{\oalign{$C$\cr \noalign{\vskip0.3ex}\hidewidth$\scriptstyle 5$\hidewidth}{}}
\def\Cf{\oalign{$C$\cr \noalign{\vskip0.3ex}\hidewidth$\scriptstyle 6$\hidewidth}{}}
\def\Cg{\oalign{$C$\cr \noalign{\vskip0.3ex}\hidewidth$\scriptstyle 7$\hidewidth}{}}

\def\Ci{\oalign{$C$\cr \noalign{\vskip0.3ex}\hidewidth$\scriptstyle 9$\hidewidth}{}}

\def\CN{\oalign{$C$\cr \noalign{\vskip0.3ex}\hidewidth$\scriptstyle N$\hidewidth}{}}
\def\CNm{\oalign{$C$\cr \noalign{\vskip0.3ex}\hidewidth$\scriptstyle N-1$\hidewidth}{}}
\def\CNmm{\oalign{$C$\cr \noalign{\vskip0.3ex}\hidewidth$\scriptstyle N-2$\hidewidth}{}}
\def\CNmmm{\oalign{$C$\cr \noalign{\vskip0.3ex}\hidewidth$\scriptstyle N-3$\hidewidth}{}}
\def\CNmmmm{\oalign{$C$\cr \noalign{\vskip0.3ex}\hidewidth$\scriptstyle N-4$\hidewidth}{}}
\def\CNmmmmm{\oalign{$C$\cr \noalign{\vskip0.3ex}\hidewidth$\scriptstyle N-5$\hidewidth}{}}

\def\Cq{\oalign{$C$\cr \noalign{\vskip0.3ex}\hidewidth$\scriptstyle q$\hidewidth}{}}
\def\Cpq{\oalign{$C$\cr \noalign{\vskip0.3ex}\hidewidth$\scriptstyle p+q$\hidewidth}{}}


{\bf Necessary and sufficient conditions for $n$-dimensional conformal Einstein spaces via dimensionally dependent identities.}

\centerline{\bf S. Brian
Edgar. }

\centerline{Department of Mathematics, }

\centerline{Link\"{o}pings universitet,}

\centerline{Link\"{o}ping,}

\centerline{Sweden S-581 83.}
\smallskip
\centerline{ email: bredg@mai.liu.se} 

\

\


{\bf Abstract.}

Listing has recently extended results of Kozameh, Newman and Tod for four-dimensional spacetimes and presented a set of necessary and sufficient conditions  for a metric to be locally conformally equivalent to an Einstein metric in  all  semi-Riemannian spaces of dimension $n\ge 4$ --- subject to a non-degeneracy restriction on the Weyl tensor.  By exploiting dimensionally dependent identities  we  demonstrate how to construct two alternative versions of these necessary and sufficient conditions  which we believe will be useful in applications. The four-dimensional case is discussed in detail and examples  are also given in five and  six dimensions. 

\

\

{\bf 1. Introduction.}

Kozameh, Newman and Tod [1] have shown that a certain pair of necessary conditions  are also sufficient for a {\it four-dimensional spacetime} to be  conformal to an Einstein spacetime --- with the exception of those spacetimes whose complex Weyl scalar invariant  $J= 0$\footnote{${}^{\dag}$}{In the conformal Einstein space with metric $g_{ab}$: $R^{ab}{}_{cd}$ represents the Riemann tensor; $C^{ab}{}_{cd}$ represents the Weyl conformal tensor; $R^{a}{}_{c}=R^{ab}{}_{cb}$ represents the Ricci tensor, $R= R^a{}_a$ the Ricci scalar,  and $\tilde R_{ab}= R_{ab}-Rg_{ab}/n$ the trace-free Ricci tensor where $2\nabla_{[a}\nabla_{b]} V^c = R_{abcd}V^d$ for an arbitrary vector $V^a$. $I,J$ are the usual complex Weyl scalar invariants in four dimensions, and $\Cb\equiv C^{ab}{}_{cd}C^{cd}{}_{ab}$. More details of the notation is given in the next section.}.  Implicit in their paper was another result: a {\it four-dimensional spacetime} with metric $g_{ab}$ can be transformed into an Einstein space  by a conformal transformation if and only if the vector $K^a$ given by 
$$K^a = 8 C^{ac}{}_{ij}\nabla^kC^{ij}{}_{ck}  \big/\ \Cb \eqno(1a)
$$
satisfies (the $n=4$ dimension version of)
$$
\tilde R_{ab} + (n-2)\Bigr(\nabla_aK_b-K_a K_b -\bigr(\nabla^cK_c - K^cK_c\bigr)g_{ab}/n\Bigr) =0
\eqno(2)$$
for the class of spacetimes where $\Cb \ne 0$. The essential ideas in [1] were  to exploit the properties that  Einstein spaces are a subset of C-spaces (spaces whose Weyl tensor is divergence-free, $\nabla_aC^{ab}{}_{cd}=0$) and that spaces conformal to C-spaces satisfy (the $n=4$ dimension version of) 
$$2 \nabla^k C^{ab}{}_{ck}+(n-3) C^{ab}{}_{ck}K^k=0
\eqno(3)$$
and hence to extract the explicit expression (1a) for  $K^a$  
by using the {\it four-dimensional} dimensionally dependent identity 
$$C^{cj}{}_{ab}C^{ab}{}_{ck} =  \delta^j_k \ \Cb/4 \ .
\eqno(4a)$$
The  principles underlying the  techniques used in  [1]  originated in a study of conformal transformations by spinor methods by Szekeres [13], and although most of the work in [1] was tensor based, some lemmas were proven by spinor methods; this is probably why Kozameh, Newman and Tod
 have commented that their method in [1] does not seem to be extendible to higher dimensions $n> 4$.  However, recently Listing [2] has exploited the same principles more generally and  shown that the condition that a particular vector 
 ${\bf T}^a$ obtained from (3) satisfies (2)  is  a necessary and sufficient condition for  conformal Einstein spaces  in {\it all semi-Riemannian spaces with dimensions $n\ge 4$} --- subject to a non-degenerate determinant condition on a matrix representation of the Weyl tensor. On the other hand, when the four dimensional results from [2] are compared with [1], the explicit form of the vector   ${\bf T}^a$  in [2] differs considerably from its counterpart $K^a$ used in [1] and  quoted above in (1a); furthermore,  the  non-degenerate determinant condition in [2] is not easily translated into a condition on the real Weyl invariant scalars in the $n$-dimensional case, although Listing states that this  condition is equivalent in four dimensions to  the complex scalar invariant Weyl scalar $J\ne0$.

It is well known that we can  write ${\bf C}^{A}{}_{C}\equiv C^{ab}{}_{cd}$ where $A\equiv[ab], C\equiv[cd]$ so that $A,C = 1,2,... , \  N(=n(n-1)/2)$
 and so consider ${\bf C}$ as an $N \times N$  trace-free matrix. (Note that  we are not making use of the symmetries, $C_{abcd}= C_{cdab}$ and $C_{a[bcd]}=0$, and in fact this construction is valid for any  trace-free double 2-form; also,  we are not defining a metric for the $N$-dimensional space, nor even making use of the $n$-dimensional space metric $g_{ab}$ in this construction.)

Instead of exploiting higher dimensional counterparts of identities such as $(4a)$, Listing's result and proof [2]  assumed $\det ({\bf C})\ne 0$ and then used the inverse matrix ${\bf C}^{-1}$ to solve (3)
for $K^a$;  however  ${\bf C}^{-1}$  cannot be easily interpreted in tensor notation without some translation and, in that form,  does not seem to be very useful in practical applications. In the illustrative example in [2] Listing restricted himself to {\it four dimensional  Riemann  space}, where he followed the  technique in [1] of using the {\it four-dimensional}  identity (4a)
to extract $K^a$, and also used properties dependent on the {\it positive definite metric}; this avoided having to deal with ${\bf C}^{-1}$ directly.

\medskip

 Although the four-dimensional identity (4a) is well known, the existence of higher dimensional analogues [3] seems less well known and one purpose of this paper is to draw attention again to the power and usefullness of such  dimensionally dependent tensor identities  [3]. (See also [4], [5], [6]). 

We shall  show in this paper how to exploit dimensionally dependent identities [3] to obtain results  valid  in  all dimensions $n\ge 4$ and automatically in all signatures; in particular we will: 

(i) reformulate Listing's results in [2], and in particular, $\det ({\bf C})\ne 0$ and the inverse matrix ${\bf C}^{-1}$ in the tensor notation of [1]; 

(ii) obtain explicit solutions for $K^a$ which avoid the  use of ${\bf C}^{-1}$ altogether.

We shall also show explicitly how the four-dimensional results implicit in  [1] can be seen as special cases of this formulation of Listing's result, and are valid for all signatures.

To demonstrate the usefulness of our versions, we will  consider the four- five -and six-dimensional cases, independent of signature.  The higher dimensional analogues [3] of the four-dimensional identity (4a)  will be the basis for our applications in five and six dimensions.  \

\

\

{\bf 2 Notation.}

We begin with some notation which we will use to  link algebraic and tensor notation, and prove a simple lemma.  Let us define, for $p \ge 1$ Weyl tensors,  a {\it chain of the zeroth kind} [6]
$$\Cp^{ab}{}_{cd}\equiv C^{ab}{}_{i_1j_1}C^{i_1j_1}{}_{i_2j_2} C^{i_2j_2}{}_{i_3j_3}  \ldots C^{i_{p-2}j_{p-2}}{}_{i_{p-1}j_{p-1}} C^{i_{p-1}j_{p-1}}{}_{cd}
$$
 noting that ${\Ca}^{ab}{}_{cd}\equiv C^{ab}{}_{cd}$.
 
A useful relation is  $$\Cp^{ab}{}_{cd}\Cq^{cd}{}_{ef}=\  \Cpq^{ab}{}_{ef} .$$
 and the following scalar invariants arise naturally,
$$\Cp \equiv \Cp^{ab}{}_{ab}.
$$
Of course there are other Weyl scalar invariants (e.g. $C^a{}_b{}^c{}_dC^b{}_e{}^d{}_f C^e{}_a{}^f{}_c$) which do not fall into this pattern.
 The simple obvious identifications which we will exploit are 
 $$\eqalign{{\bf C}^p & = \Cp^{ab}{}_{cd}\cr 
 <\!\! {\bf C}^p\!\!> & = \Cp}\eqno(5)$$ 
where $<\ >$ denotes the trace of a matrix; but note that expressions like $\Cb^{cj}{}_{ci}$ have no such obvious identification in the matrix notation.
 
 The Cayley-Hamilton theorem for the trace-free $N\times N$ matrix $ {\bf C}$ is given by,
$$c_0{\bf C}^N + c_2{\bf C}^{N-2}+  c_3{\bf C}^{N-3}+\ldots  + c_{N-2}{\bf C}^2+ c_{N-1}{\bf C}+  c_N {\bf I}=0
\eqno(6)$$
where ${\bf I}$ is the $N \times N$ identity matrix, and 
$$\eqalign{  c_0 & =  1,  \  c_2 = -{1\over 2}<\!\!{\bf C}^2\!\!>, \  c_3=  -{1\over {3}}<\!\!{\bf C}^3\!\!>, \  c_4= -{1\over {4}}\Bigl(<\!\!{\bf C}^4\!\!> - {1\over 2}<\!\!{\bf C}^2\!\!>^2\Bigl), \cr  c_5 & = - {1\over {5}}\Bigl(<\!\!{\bf C}^5\!\!> - {5\over 6}<\!\!{\bf C}^2\!\!><\!\!{\bf C}^3\!\!>\Bigl),  \cr  c_6 &=  - {1\over 6}\Bigl(<\!\! {\bf C}^6\!\!> - {3\over 4} <\!\! {\bf C}^2\!\!><\!\! {\bf C}^4\!\!>
-{1\over 3} <\!\! {\bf C}^3\!\!>^2 +{1\over 8} <\!\! {\bf C}^2\!\!>^3\Bigr), 
\cr
& \ldots , \cr
c_N & = -{1\over N}\Bigl( <\!\! {\bf C}^N\!\!> +  \ldots \qquad +  \ldots \qquad \Bigr) }\eqno(7)$$ 
are the usual characteristic coefficients ; since ${\bf C}$ is trace-free, $c_1=0$.
This theorem can easily be rewritten in chain notation as
$$   c_0 \CN^{ab}{}_{cd}+ c_2 \  \CNmm^{ab}{}_{cd} +c_3 \  \CNmmm^{ab}{}_{cd}+ \ldots  + c_{N-2} \  \Cb^{ab}{}_{cd}  + c_{N-1} \ \Ca^{ab}{}_{cd} + c_{N} \delta^{a}_{[c}\delta^{b}_{d]}=0
 \eqno(6') $$ 
where the characteristic coefficients  are now given in terms of Weyl scalar invariants  by 
$$\eqalign{  c_0 & =  1,  \   \   c_2 = -{1\over 2}\Cb, \    \  c_3=  -{1\over {3}} \Cc, \     c_4= -{1\over {4}}\Bigl( \Cd - {1\over 2} \Cb^2\Bigl), \cr  c_5 & =  -{1\over {5}}\Bigl( \Ce - {5\over 6} \Cb \Cc\Bigl),  \cr  c_6 &=  - {1\over 6}\Bigl(\Cf- {3\over 4}  \Cb \Cd
-{1\over 3} \Cc^2 +{1\over 8} \Cb^3\Bigr), \cr
& \ldots , \cr
c_N & =- {1\over N}\Bigl( \CN + \ldots +  \qquad \ldots \qquad \Bigr ) 
}\eqno(7')$$
From  the well-known result  
$(-1)^N \det({\bf C}) = c_N
$
the required translation of $\det ({\bf C})\ne 0$ into tensor language follows immediately,
$$ 0 \ne (-1)^{N+1} N \det ({\bf C})=  \CN + \ldots \qquad + \ldots
\ \eqno(8)$$

\smallskip
{\bf Lemma.}   In $n$-dimensional spaces, the inhomogeneous algebraic equation for the vector $V^d$
$$C^{ab}{}_{cd}V^{d}=H^{ab}{}_{c}
$$
has  a unique solution when 
 condition (8) holds; 
the solution is given by
$$   V^a =  {2\over (n-1) c_{N}} H^{ij}{}_b \Bigl(c_0 \  \CNm^{ab}{}_{ij}+ c_2\   \CNmmm^{ab}{}_{ij}+c_3 \  \CNmmmm^{ab}{}_{ij}+ \ldots  + c_{N-2} \ C^{ab}{}_{ij} \Bigr)
  $$ 
where $N=n(n-1)/2$.
The coefficients $c_0,c_2,c_3, \  \ldots c_{N-2}, c_N$ are the usual characteristic coefficients of the Cayley-Hamilton Theorem given in $(7')$. \bigskip

{\bf Proof.}
We consider the Cayley-Hamilton theorem for the $N \times N$  trace-free matrix ${\bf C}$
 in tensor notation in $(6')$, with characteristic coefficients given by $(7')$.

Multiplying by $V_a$   gives
$$ \eqalign{  0  & = V_a\Bigl(  c_0 \  \CN^{ab}{}_{cd}+ c_2\   \CNmm^{ab}{}_{cd}+c_3 \  \CNmmm^{ab}{}_{cd}+c_4  \  \CNmmmm^{ab}{}_{cd}+   \cr & \qquad \qquad\qquad \qquad \ldots  + c_{N-2}  \  \Cb^{ab}{}_{cd}  + c_{N-1}  \  \Ca^{ab}{}_{cd}\Bigr) + c_{N} V_{[c}\delta^{b}_{d]}
\cr & 
 = V_a \Ca^{ab}{}_{ij}\Bigl(  c_0  \  \CNm^{ij}{}_{cd}+ c_2  \  \CNmmm^{ij}{}_{cd}+c_3  \  \CNmmmm^{ij}{}_{cd}+c_4  \  \CNmmmmm^{ij}{}_{cd}+ \cr & \qquad \qquad\qquad \qquad  \ldots  + c_{N-2}  \  \Ca^{ij}{}_{cd} \Bigr) + c_{N-1}  \  V_a \Ca^{ab}{}_{cd} + c_{N}  \  V_{[c}\delta^{b}_{d]}
 }$$

From which, by taking the trace and remembering $\Ca^{ab}{}_{cd}\equiv C^{ab}{}_{cd}$ is trace-free, we obtain, 
$$ \eqalign{  0  & =  V_a C^{ab}{}_{ij}\Bigl(  c_0  \  \CNm^{ij}{}_{cb}+ c_2  \  \CNmmm^{ij}{}_{cb}+c_3  \  \CNmmmm^{ij}{}_{cb}+c_4  \  \CNmmmmm^{ij}{}_{cb}+ \cr & \qquad \qquad\qquad \qquad  \ldots  + c_{N-2}  \  C^{ij}{}_{cb} \Bigr)  +{n-1\over 2} c_{N} V_c
 }$$

Rearranging gives the solution in the lemma. 

\hfill \rightline{$\square$}

\medskip

 For future reference, we note that the {\it four-dimensional identity}  $(4a)$ can be written in the chain notation as
 $$\Cb^{cj}{}_{ck}= \delta^j_k \Cb/4
$$
and this is actually a special case of the more general identity [6] 
 {\it in  four dimensions only}
$$\Cp^{cj}{}_{ck}= \delta^j_k \Cp/4.  \qquad p=2,3,4 \ \ldots  
\eqno(4^*)$$
We have preferred the notation $\Cb, \Cc, \ldots \ \Cp, \ldots \ $ to the (possibly confusing) notation $C^2, C^3, \ldots \ C^p, \ldots \ $ used in [1], [2] and elsewhere for these Weyl scalar invariants.

\medskip

Finally we note a very useful dimensionally {\it independent} identity (a direct consequence of the first Bianchi identity),
$$4C_{a[ij]b}C^{cijd}= C_{abij}C^{cdij}
$$
 
\

\

{\bf 3. Reformulating Listing's result and rederiving the implicit results in [1].}

With this Lemma we determine the vector $K^a$ from (3),
$$   K^a =   {4\over (n-3)(n-1)c_N} \nabla^kC^{ij}{}_{kb} \Bigl(c_0 \   \CNm^{ab}{}_{ij}+ c_2  \  \CNmmm^{ab}{}_{ij}+c_3  \  \CNmmmm^{ab}{}_{ij}+ \ldots  + c_{N-2}  \  C^{ab}{}_{ij} \Bigr)
  \eqno(9)$$ 
providing restriction (8) holds, where $N=n(n-1)/2$. 

Substituting this value for $K^a$ back into (2)  gives necessary and sufficient conditions defined only in terms of the geometry.

We can therefore reformulate Theorem 4.5 in [2] as follows:
\medskip
{\bf Theorem 1.}  {\it A semi-Riemannian manifold, with a  Weyl tensor subject to the restriction $(8)$, is locally conformally related to an Einstein space if and only if the vector field $K^a$ given in (9) satisfies (2)}.

\medskip
\underbar{Four dimensions.}

We now will retrieve the four-dimensional results implicit in [1]   from this theorem in a signature independent manner. 
Instead of substituting for $K^a$ with the expression (1a) as used in [1],  or the four-dimensional version of the expression involving ${\bf C}^{-1}$ in [2], we can now   use the $n=4$ version of (9),
$$   K^a =   {4\over 3c_6} \nabla^kC^{ij}{}_{kb} \Bigl(c_0 \Ce^{ab}{}_{ij}+ c_2 \Cc^{ab}{}_{ij}+c_3 \Cb^{ab}{}_{ij}+  c_{4} C^{ab}{}_{ij} \Bigr)
  \eqno(10)$$ 
providing $c_6\ne 0$ with the characteristic coefficients  given by $(7')$ in terms of Weyl scalars. 

However, since the solution for $K^a$ is unique we should be able to see precisely the links between these two expressions in (1a) and (10).
To retrieve the result (1a) from (10) we simply substitute (3a) back into the right hand side of (a slightly rearranged)  (10), for all terms except the last one, and use identity (4a) on each of these terms,
 $$  \eqalign {3c_6K^a  & =   -{2} K^k C^{ij}{}_{kb} \Bigl( \Ce^{ab}{}_{ij}+ c_2 \Cc^{ab}{}_{ij}+c_3 \Cb^{ab}{}_{ij} \Bigr) +{4} c_4C^{ab}{}_{ij}\nabla^kC^{ij}{}_{kb}\cr & =
-2 K^k  \Bigl( \Cf^{cb}{}_{kb}+ c_2 \Cd^{cb}{}_{kb}+c_3 \Cc^{cb}{}_{kb} \Bigr) +4 c_4C^{ab}{}_{ij}\nabla^kC^{ij}{}_{kb} \cr &  
= -{1\over 2} K^a  \Bigl(\Cf+ c_2 \Cd+c_3 \Cc \Bigr) +4 c_{4} C^{ij}{}_{cb} \nabla^kC^{ij}{}_{kb} 
} $$
which rearranges to
$$   (6c_6 +  \Cf+c_2 \Cd + c_3 \Cc)K_c=  8 c_{4} C^{ij}{}_{cb} \nabla_kC^{kb}{}_{ij}  $$
  and via $(7')$ to 
  $$c_4\Cb K^a=  -8 c_{4} C^{ab}{}_{ij} \nabla^kC^{ij}{}_{kb}  
  $$
  and hence to (1a) --- providing $\Cb\ne 0\ne \Cd-\Cb^2/2$. 
  
 Kozameh, Newman and Tod [1]  have shown that it is possible to obtain three alternative versions to (1a) for $K^a$ consisting of
  $$K^a   = 8 \Cb^{ac}{}_{ij}\nabla^kC^{ij}{}_{ck}  \big/ \Cc \qquad \hbox{for} \ \Cc\ne 0
  $$ 
  together with two other similar expressions with invariants involving the dual of the Weyl tensor. We do not wish to depend on  expressions with duals in this work since we wish to generalise the four-dimensional case to all higher dimensions.   Instead, from the general identity $(4^*)$ we see that we can obtain the various expressions for $K^a$ for all integers  $p\ge 2$
  $$  K^a   = 8 \ \Cpm^{ac}{}_{ij}\nabla^kC^{ij}{}_{ck}  \big/ \Cp\qquad \hbox{for} \ \Cp\ne 0 \eqno(1^*).
$$ 
 
 So, in an analogous manner to which we retrieved $(1a)$ from $(10)$, we can also retrieve $(1^*)$ for $p=3,4,6$; by use of the Cayley-Hamilton theorem we could also retrieve the results for $p=5$ and $p\ge 7$, although after $p=6$ these expressions are not independent precisely because of the Cayley-Hamilton theorem. We note that we cannot retrieve the version of $(1^*)$ directly with $p=5$  from (10) due to the fact that the relevant term  involved the coefficient $c_1$ which is identically zero, being the trace of $C_{ab}{}^{cd}$. 
 
 It is well known that in  four dimensions there exist only four algebraically independent  Weyl scalar invariants: these are usually given in terms of the well known complex  invariants $I$ and $J$, which come naturally from considerations in terms of spinors or complex Weyl tensors (which exploit the very simple and unique  structure of the Weyl dual in four dimensions).  If we wish to only employ structures which can be exploited  in higher dimensions, we would consider instead the four algebraically independent scalars $\Ca,\Cb,\Cc,\Cd$ (or any four from $\Ca,\Cb,\Cc,\Cd,\Ce$) [7] (see also [8],[9]); and providing there is at least one Weyl scalar invariant which is non-zero, we can calculate $K^a$ explicitly.
 
\medskip
\underbar{Higher dimensions.}  

Theorem 1 gives an explicit  tensor expression for $K^a$ for all dimensions. So, for example in {\it five dimensions}, $N=10$ and we obtain
 $$   K^a =   {1\over 2c_{10}} \nabla^kC^{ij}{}_{kb} \Bigl(c_0 \Ci^{ab}{}_{ij}+ c_2 \Cg^{ab}{}_{ij}+c_3 \Cf^{ab}{}_{ij}+ \ldots + c_{7} \Cb ^{ab}{}_{ij}+ c_{8} C^{ab}{}_{ij}  \Bigr)
 $$ 
providing $c_{10}\ne 0$ with the characteristic coefficients  given by $(7')$ in terms of Weyl scalars. 
Clearly  this expression involves quite high order terms in the Weyl tensor, and we note that in the four-dimensional case, we were able to get lower order expressions for $K^a$ by exploiting the individual dimensionally dependent identities $(4^*)$ rather than the Cayley-Hamilton Theorem implicit in Theorem 1.
So,  in higher dimensions, we would also expect to exploit  individual dimensional identities, analogous to the four-dimensional identities  $(4^*)$, 
in order to obtain alternative lower order forms in the Weyl tensor for $K^a$ in $n$ dimensions. 

\

\

{\bf 4. Lower order versions of $K^a$ in $n$-dimensions.}

In four dimensions, providing we pay the price that at least one scalar invariant is non-zero, we are able to solve the four-dimensional version of $(3)$, and obtain  $K^a$. The versions from the individual  identities $(4^*)$ are  more concise and manageable than  the  versions from the Cayley-Hamilton theorem; on the otherhand, in four dimensions, the latter include all the simpler results  as special cases, and are more general, in the sense of  weaker restrictions on the Weyl scalar invariants.  Of course also  the Cayley-Hamilton approach is applicable in all dimensions. However, we now  will show how to obtain alternative simpler versions  for $K^a$ in dimensions other than four, by direct application of dimensionally dependent identities analogous to those used in four dimensions in [1].

The set of identities $(4^*)$ are all consequences of the well-known four-dimensional identity\footnote{${}^{\dag}$}{It is interesting to note, in this notation, that the underlying identity for the Cayley-Hamilton theorem for the trace-free $6\times 6$ matrix ${\bf C}$ is ${\bf C}^{[A}{}_{[P}\delta^{BCDEF]}_{QRSTU]}=0$ [3].} [5]
 $$C^{[ab}{}_{[cd}\delta^{e]}_{f]}=0\eqno(11)$$
and we shall now exploit their higher dimensional analogues in the same way as in [1] for the four-dimensional case:  for five-dimensions 
 we have $$C^{[ab}{}_{[cd}\delta^{ef]}_{hi]}=0 \eqno(12)$$
and for six-dimensions
 $$C^{[ab}{}_{[cd}\delta^{efg]}_{hij]}=0  \eqno(13)$$
and so on [3].

These lead to 
$$C^{cd}{}_{pe}C^{ph}{}_{ab}C^{[ab}{}_{[cd}\delta^{ef]}_{hi]}=0 \eqno(14)$$
for five dimensions, and 
$$C^{cd}{}_{ef}C^{hi}{}_{ab}C^{[ab}{}_{[cd}\delta^{efg]}_{hij]}=0 \eqno(15)$$
for six dimensions, and so on. Of course, these are just representative for each dimension of the various identities that can be created; but it is of interest to note that for the lowest order in Weyl in each dimension --- in five dimensions cubic, in six dimensions cubic, in seven dimensions quartic,and so on ---  there is only one possibility, but as we look at higher orders --- in five dimensions quartic, in six dimensions quartic, in seven dimensions quintic, and so on ---  there will be many more possibilities.
Unfortunately, when expanded, such   two-index identities in higher dimensions  will not have such a simple form as  the four-dimensional identities $(4^*)$. 

In   higher dimensions there exist  greater numbers of  algebraically independent Weyl scalar invariants, and so there will be  greater numbers of two-index tensor identities analogous to $(4^*)$ for each dimension $n> 4$;  furthermore,  the   tensor identities analogous to $(4^*)$ will be based on (14), (15),  ...  and so will also require more Weyl tensors as the dimension increases. Since, in higher dimensions,  a product of three or more Weyl tensors yields more than one Weyl scalar algebraically independent  of each other and of invariants of lower order (e.g. in general, $\Cc$ and  $C^a{}_i{}^b{}_j C^i{}_c{}^j{}_d C^c{}_a{}^d{}_b $ are algebraically independent in $n\ge 6$ dimensions), and more than one  algebraically independent two-index tensor, we expect the higher dimensional analogues of identities  $(4^*)$ to consist of linear combinations of Weyl tensors on both sides of the identity. 

However, it is important to note that, although these higher order two-index identities will have more terms, and higher products,  they  will have the same crucial structure, which we can represent by
$$L\{\Cp\}^j{}_k= \delta^j_k L\{\Cp\}/n
\eqno(16)$$
where $L\{\Cp\}^j{}_k$  represents  a two-index tensor consisting of a linear combination of products of $p$ Weyl tensors, and  $L\{\Cp\}\equiv L\{\Cp\}^i{}_i$ represents  a linear combination of scalar products of $p$ Weyl tensors. (Note that each term is unlikely to be just a simple chain.) It follows, for $L\{\Cp\}\ne 0$ from (16)
$$K^j= n L\{\Cp\}^j{}_kK^k\big/L\{\Cp\}
\eqno(17)$$
and hence all the terms involving the vector $K^k$ on the right hand side --- which will each contain a factor of the form $C^{. .}{}_{. k}K^k$ --- can be replaced using (3) by 
$$C^{. .}{}_{. k}K^k= - {2\over n-3} \nabla^kC^{. .}{}_{. k}
\eqno(18)$$

\

We can summarise these new results as follows,
\medskip
{\bf Theorem 2. }  {\it An $n$-dimensional semi-Riemannian manifold    with  a nondegenerate Weyl tensor restricted by $L\{\Cp\} \ne 0$ (where $L\{\Cp\} $ is associated with an identity of the form (16))    is locally conformally related to an Einstein space if and only if the vector field $K^a$ given in (17), with the appropriate substitutions (18), satisfies (2)}.

\medskip

Clearly from Theorem 2 we will obtain for $K^a$ much lower order expressions in the Weyl tensor than from Theorem 1; for example, in five dimensions using Theorem 1 will require terms involving products of ten Weyl tensors, whereas if we use Theorem 2 it looks possible to use terms with only three Weyl tensors, from (14).

\

\

 {\bf 5. Five- and six-dimensional spaces.}

For higher dimensional spaces we can use Theorem 1 with the respective substitutions $n=5,6,\  \ldots $ into (9).  But for spaces where we know dimensionally dependent identities of the form $(14), (15) \ldots$ we can use Theorem 2.

So for spaces with dimension $n>4$ there is the need to systematically write out explicitly the two-index identities such as $(14), (15) \ldots$  for $n=5,6, \ldots \ $.  In this section we give just a few examples in five and six dimensions as illustrations.

\underbar{Five dimensions.}

After a straightforward  calculation we find that $(14)$ expands as follows,
$$\eqalign{C^{aj}{}_{bc}C^{bc}{}_{de}C^{de}{}_{ak} & - 2 C^{aj}{}_{bk}C^{bc}{}_{de}C^{de}{}_{ac} -4 C^{aj}{}_{bc}C^{bd}{}_{ek}C^{ce}{}_{ad}\cr
 & ={1\over 5}\bigl(C^{ab}{}_{cd}C^{cd}{}_{ef}C^{ef}{}_{ab}-4C^{ab}{}_{cd}C^{ce}{}_{af}C^{df}{}_{be}\bigr) \delta^j_k
}\eqno(19)$$
Although this appears to have the structure of $(16)$, unfortunately when we also consider the scalar identity closely related to $(14)$ 
$$C^{cd}{}_{ab}C^{ki}{}_{ej}C^{[ab}{}_{[cd}\delta^{ej]}_{ki]}=0 $$
we find 
$$C^{ab}{}_{cd}C^{cd}{}_{ef}C^{ef}{}_{ab}=4C^{ab}{}_{cd}C^{ce}{}_{af}C^{df}{}_{be}
\eqno(20)$$
(This five-dimensional scalar identity was also noted in [14], where it was obtained from the {\it five dimensional} identity $C^{ab}{}_{[cd}C^{cd}{}_{ef}C^{ef}{}_{ab]}\equiv 0$.)
This means that $(19)$ does not have the structure of $(16)$ as we hoped, since its right hand side is identically zero. However, we do have an interesting two-index {\it five-dimensional} identity which will be useful in other contexts,
$$\eqalign{C^{aj}{}_{bc}C^{bc}{}_{de}C^{de}{}_{ak}  - 2 C^{aj}{}_{bk}C^{bc}{}_{de}C^{de}{}_{ac} -4 C^{aj}{}_{bc}C^{bd}{}_{ek}C^{ce}{}_{ad} =0
}\eqno(21)$$

On the otherhand,  if we consider one (of a number of) quartic identity in five dimensions,
$$C^{qg}{}_{ip}C^{ip}{}_{qe}C^{ab}{}_{cd}C^{[cd}{}_{[ab}\delta^{ef]}_{gh]}=0
\eqno(22)$$
we obtain 
$$\eqalign{5C^{qj}{}_{ip}C^{ip}{}_{qk}& C^{ab}{}_{cd}C^{cd}{}_{ab} -  8 C^{qg}{}_{ip}C^{ip}{}_{qk}C^{ab}{}_{cg}C^{cj}{}_{ab}
+ 8 C^{qg}{}_{ip}C^{ip}{}_{qe}C^{be}{}_{ag}C^{aj}{}_{bk}\cr &
-4 C^{qg}{}_{ip}C^{ip}{}_{qe}C^{ab}{}_{gk}C^{ej}{}_{ab} 
-8 C^{qg}{}_{ip}C^{ip}{}_{qe}C^{ae}{}_{bk}C^{bj}{}_{ag}
\cr & \qquad \qquad \qquad =\bigl(C^{qg}{}_{ip}C^{ip}{}_{qg}C^{ab}{}_{cd}C^{cd}{}_{ab} -4C^{qg}{}_{ip}C^{ip}{}_{qd}C^{ab}{}_{cg}C^{cd}{}_{ab}\bigr)\delta^j_k
}\eqno(23)$$
which again appears to have the structure of  (16). (The use of the identity (21) does not give any simplification.)  We need to determine whether the right hand side of (23) is non-zero.  Unlike in the cubic case where there was only one possible scalar identity $(20)$, there will be a number of quartic scalar identities in five dimensions ; so although we will return to examine these another time, for now we simply note, via a  counterexample, that the right hand side  cannot be identically zero [10].

Hence, via the substitution (18),  we obtain the following form for $K^j$,
$$\eqalign{K^j & = \cr & \   -\Bigl( 5C^{qj}{}_{ip} C^{ab}{}_{cd}C^{cd}{}_{ab}\nabla^k C^{ip}{}_{qk}  -  8 C^{qg}{}_{ip}C^{ab}{}_{cg}C^{cj}{}_{ab}\nabla^k C^{ip}{}_{qk}
+ 8 C^{qg}{}_{ip}C^{ip}{}_{qe}C^{be}{}_{ag} \nabla^k C^{aj}{}_{bk}\cr &
\qquad \qquad\qquad \qquad -4 C^{qg}{}_{ip}C^{ip}{}_{qe}C^{ej}{}_{ab} \nabla^k C^{ab}{}_{gk}
-8 C^{qg}{}_{ip}C^{ip}{}_{qe}C^{bj}{}_{ag}\nabla^k C^{ae}{}_{bk}\Bigr)\cr &  \qquad \qquad\qquad\qquad \qquad\qquad \Big/ \bigl(C^{qg}{}_{ip}C^{ip}{}_{qg}C^{ab}{}_{cd}C^{cd}{}_{ab} -4C^{qg}{}_{ip}C^{ip}{}_{qe}C^{ab}{}_{cg}C^{ce}{}_{ab}\bigr)
}\eqno(24)$$
 We can then substitute this value for $K^j$ into the five dimension version of (2),
$$
\tilde R_{ab} + 3\Bigr(\nabla_aK_b-K_a K_b -\bigr(\nabla^cK_c - K^cK_c\bigr)g_{ab}/5\Bigr) =0
\eqno(2b)$$
to obtain the required necessary and sufficient condition.

 \medskip
\underbar{Six dimensions.}

When we expand $(15)$  we obtain,
$$C_{ak}{}^{bc}C^{aj}{}_{de}C_{bc}{}^{de} - 2 C_{ak}{}^{bj}C^{ac}{}_{de}C_{bc}{}^{de} - 4C_{ak}{}^{bc}C^{dj}{}_{be}C_{cd}{}^{ae}= {1\over 6}\bigl(C_{ab}{}^{cd}C_{cd}{}^{ef}C_{ef}{}^{ab}-4C_{ab}{}^{cd}C^{ae}{}_{cf}C_{de}{}^{bf}\bigr)\delta_k^j
\eqno(25)$$
Unlike for the case of the five-dimensional two-index identity (19), there is no related scalar identity  cubic in Weyl tensors,  analogous to  (20).  So therefore we do not need to worry about the possibility of the right hand side of (25) being zero, and so  we have an identity which has  precisely the structure  (16).

Hence, via the substitution (18),  we obtain the following form for $K^j$,
$$\eqalign{K^j & = -4 \bigl(C^{aj}{}_{de}C_{bc}{}^{de}\nabla^k C_{ak}{}^{bc} - 2 C^{ac}{}_{de}C_{bc}{}^{de}\nabla^k C_{ak}{}^{bj} - 4C^{dj}{}_{be}C_{cd}{}^{ae}\nabla^k C_{ak}{}^{bc}\bigl) \cr & \qquad\qquad  \big/\bigl(C_{ab}{}^{cd}C_{cd}{}^{ef}C_{ef}{}^{ab}-4C_{ab}{}^{cd}C^{ae}{}_{cf}C_{de}{}^{bf}\bigr)}
\eqno(26)$$
 We can then substitute this value for $K^j$ into the six-dimensional version of (2),
$$
\tilde R_{ab} + 4\Bigr(\nabla_aK_b-K_a K_b -\bigr(\nabla^cK_c - K^cK_c\bigr)g_{ab}/6\Bigr) =0
\eqno(2c)$$
to obtain the required necessary and sufficient condition.

\medskip
There do not seem to be many explicit examples of identities for the Weyl tensor in higher dimensions in the literature. However, there does exist a six-dimensional two-index tensor identity quartic in the Weyl tensor which was identified some time ago by Lovelock [4]. 
A double three-form with the antisymmetric and trace-free properties respectively
$$\eqalign{H_{abk}{}^{def} & = H_{[abk]}{}^{[def]}\cr
H_{abi}{}^{dei} & = 0}
\eqno(27)$$
in six dimensions (and lower) satisfies the identity
$$H_{abk}{}^{def}H_{def}{}^{abj}={1\over 6}\delta^j_k H_{abc}{}^{def}H_{def}{}^{abc}\ .
\eqno(28)$$
 This then becomes a quartic identity for Weyl with the choice
 $$\eqalign{H_{ijk}{}^{abc}  =A_{ijk}{}^{abc}- {9\over 8}A_{r[jk}{}^{r[bc} \delta^{a]}_{i]} + {3}A_{rs[i}{}^{rs[c} \delta^{a}_{i} \delta^{b]}_{k]}  
- {1\over 4}A_{rst}{}^{rst} \delta^{[a}_{[i} \delta^{b}_{j}\delta^{c]}_{k]}  }\eqno(29)$$
where
$$A_{ijk}{}^{abc}=4C_{[ij}{}^{h[a}C_{k]h}{}^{bc]}
\eqno(30)$$

By substituting (29) into (28) we obtain 
$$ \eqalign{A_{abk}{}^{cde}A_{cde}{}^{abj} & + 3A_{abk}{}^{abc}A_{cde}{}^{dej}+6A_{abk}{}^{acd}A_{cde}{}^{bej}-3A_{abc}{}^{ade}A_{dek}{}^{bc j} - A_{abc}{}^{abc}A_{dek}{}^{dej}+6A_{abc}{}^{abd}A_{dek}{}^{cej}\cr & = \delta^j_k\Bigl(A_{abc}{}^{abc}A_{def}{}^{def}+9A_{abc}{}^{ade}A_{def}{}^{bcf}-9A_{abc}{}^{abd}A_{def}{}^{cef}
-A_{abc}{}^{def}A_{def}{}^{abc}\Bigr)/6}\eqno(31)$$
which is precisely the structure of (16). Hence,   we obtain the following form for $K^j$,
$$ \eqalign{K^j & = \cr  & \  -4 \Bigl(K^k A_{kab}{}^{cde}A_{cde}{}^{abj}  + 3K^kA_{kab}{}^{abc}A_{cde}{}^{dej}+6K^kA_{kab}{}^{acd}A_{cde}{}^{bej}-3A_{abc}{}^{ade}K^kA_{kde}{}^{bc j}\cr & \qquad \qquad\qquad - A_{abc}{}^{abc}K^kA_{kde}{}^{dej}+6A_{abc}{}^{abd}K^kA_{kde}{}^{cej}\Bigr)\cr & \qquad\qquad \qquad\qquad  \big / \bigl (A_{abc}{}^{abc}A_{def}{}^{def}+9A_{abc}{}^{ade}A_{def}{}^{bcf}-9A_{abc}{}^{abd}A_{def}{}^{cef}
-A_{abc}{}^{def}A_{def}{}^{abc}\bigr)}\eqno(32)$$
where all terms involving $K^k$ on the right hand side are replaced via the substitution
$$K^k A_{kab}{}^{cde}= 2\nabla^kC_{k[a}{}^{p[c}C_{b]p}{}^{de]}- \nabla^kC_{kp}{}^{[de}C_{ab}{}^{c]p}
\eqno(33)$$
and all other terms replaced with (30).
 We then substitute the value for $K^j$ from (32) into $(2c)$ to obtain the required necessary and sufficient condition. It is of course necessary to check that the right hand side of $(32)$ is not identically zero; this can be confirmed with a simple counterexample [10].

\

\

{\bf 6. Summary and Discussion.} 

We have demonstrated the power of  dimensionally dependent identities (either as the Cayley-Hamilton Theorem, or as individual identities in specific dimensions such as (4*), (23), (25), (31)) to translate the existence result of Listing [2] into  versions which can be applied directly; furthermore, these applications are completely independent of signature. Theorem 1 gives a general reformulation of Listing's result on the existence of necessary and sufficient conditions in a form  which can be directly exploited.  So if one wishes to test if a particular metric is conformally Einstein, then it is a simple procedure to find its Weyl and trace-free Ricci tensor and to test directly --- for instance using GRTENSOR II [10] --- if condition (2) is satisfied.

   We have drawn attention to the higher dimensional analogues $ (23), (25), (31), \ldots \ $ of identities (4*), which are the basis for Theorem 2; for applications, Theorem 2 would seem to provide a simpler and more manageable tool --- providing the appropriate identities are known.

A major complication when we move to  dimensions $n>4$ is that there are many more Weyl invariant scalars, and of course they cannot all be written in the form $\Cp$. The same sort of detailed analysis of the Weyl  invariant scalars for $n=5,6, \ldots$ as has been (partly) carried out for $n=4$, as well as a systematic presentation of two-index identities, is a necessary prerequisite for a systematic examination of all possible versions of the vector $K^a$.    As Bonanos has pointed out  [11] existing detailed studies such as [12] do not  take into account invariants formed from duals, or identities  from the Cayley-Hamilton Theorem.

There are still a number of interesting issues to be investigated further.

There appears to be an important difference between even and odd dimensions: for even dimensions  $n=4,6,8, \ldots $ the simplest two-index identity involves products of $2,3,4,\ldots$ Weyl tensors respectively, such as (4a) and  (25), with a 'delta term' on the right hand side; for odd $n=5,7, \ldots $ the simplest two-index identity involves products of $3,4,\ldots$ Weyl tensors respectively, such as (19), but it would seem likely that as in (21) the 'delta term' on the right hand side disappears because of an identically zero coefficient. For the investigations in this paper we need an identity of the former type, so in general it appears that for even dimensions $n=2m$ we will be able to exploit comparatively simple identities involving products of $m$ Weyl tensors, while for odd dimensions $n=2m-1$ we will only have more complicated identities involving products of $m+1$ Weyl tensors.  On the otherhand, we anticipate that in other investigations the simple identities such as (21) will be very useful.

 Listing [2] has stated that the condition $\det{\bf C}\ne 0$ {\it in four dimensions} is equivalent to the complex Weyl invariant scalar $J\ne 0$; by a little manipulation this can be shown to be equivalent to at least one of $\Cc$ and $\Ce$ being non-zero. It would be useful to know this condition in higher dimensions in terms of the real Weyl invariant scalars, and hence understand it better.
The fact that the right hand side of identity (19) is identically zero should alert us to the possibility of identically zero scalars arising in some situations.

The use of the dual Weyl tensor makes work in four dimensions comparatively easy--- for instance there is a basis of four Weyl scalar invariants none of which is higher than cubic in Weyl, compared to having a basis with invariants up  to fifth order in Weyl if the dual tensor is not used; in higher dimensions the major advantage  (the dual Weyl tensor is also a double two form) does not apply, and work gets more complicated. However, we believe that it is still possible to take advantage of other benefits of the dual tensor, and we will discuss this possibility, together with the other points mentioned here,  elsewhere.

 Finally, we note that the necessary and sufficient conditions investigated by Listing [2] and in this paper were different from the necessary and sufficient conditions investigated explicitly in [1] in four dimensions. In [1] these conditions involved the Bach tensor which of course is only defined in four dimensions. It will be shown elsewhere how the techniques in [2] and in this paper can be used to investigate these alternative conditions as well as to generate an $n$-dimensional generalisation of the Bach tensor.

\

\

{\bf Acknowledgements.} Thanks to Jonas Bergmann for help with the five and six dimensional identities, and to Magnus Herberthsson and Jos\'{e} Senovilla for discussions.

\

\

{\bf References.}

1. Kozameh, C.N., Newman, E.T. and Tod, K.P.:  Conformal Einstein spaces. {\it Gen. Rel. Grav.}, {\bf 17} (1985), 343-352.

2. Listing, M.:  Conformal Einstein Spaces in $N$-Dimensions, {\it Annals of Global Analysis and Geometry}, {\bf 20} (2001), 183-197.

3. Edgar, S.B. and H\"oglund, A.:  Dimensionally dependent tensor
identities by double antisymmetrisation.  {\it  J. Math. Phys.}, {\bf 43} (2002), 659-677.

4. Lovelock, D.: Dimensionally dependent identities, {\it Proc. Camb. Phil. Soc.},{\bf 68} (1970), 345-350.

5. Edgar, S.B. and Wingbrant, O.:  Old and new results for superenergy tensors using dimensionally
dependent identities.  {\it  J. Math. Phys.}, {\bf 44} (2003), 6140 - 6159.

6. Wingbrant, O.: Highly structured tensor identities for (2,2)-forms in four dimensions. \ {\it Preprint}:  http://arxiv.org/abs/gr-qc/0310120

7. Narlikar, V.V. and  Karmarkar, K.R.:  The scalar invariants of a general relativity metric. {\it Proc. Indian Acad.  Sci.} A , {\bf 29} (1949), 91.

8. Sneddon, G. E.: On the algebraic invariants of the four dimensional Riemann tensor. {\it Class. Quantum Grav.}, {\bf 3}  (1986), 1031-1032.

9. Harvey, A.: On the algebraic invariants of the four dimensional Riemann tensor.
{\it Class. Quantum Grav.}, {\bf 7} (1990), 715-716.

10. The counterexamples are obtained by confirming the existence of non-zero  Weyl scalars for a particular metric; and we can do this very easily with GRTENSOR II,  available at  http://grtensor.org .

11. Bonanos, S.:  A new spinor identity and the vanishing of certain Riemann tensor invariants. {\it Gen. Rel. Grav.}, {\bf 30} (1998), 653-658.

12. Fulling, S.A. {\it et al}: Normal forms for tensor polynomials: I. The Riemann tensor. {\it Class. Quantum Grav.}, {\bf 9}  (1992), 1151-1197.

13. Szekeres, P.:  Spaces conformal to a class of spaces in general relativity.  \ {\it Proc. Roy. Soc.  Lon.} A , {\bf 274} (1963), 206-212.

14. Jack, I. and Parker, L.: Linear independence of renormalisation counterterms in curved space-times of arbitrary dimensionality.  {\it  J. Math. Phys.}, {\bf 28} (1987), 1137 - 1139.

\end

  from which we obtain directly its trace
  $$   c_0 \CN^{a}{}_{c}+ c_2 \CNmm^{a}{}_{c}+c_3 \CNmmm^{a}{}_{c}+c_4 \CNmmmm^{a}{}^{c}+ \ldots  + c_{N-2} \Cb^{a}{}_{c} + {N-1\over 2}  c_{N} \delta^{a}_{c}=0
  $$   
where $c_0 = 1, c_2 = -\Cb/2, c_3=  -\Cc/3, c_4= -\Bigl(\Cd-\Cb^2/2\Bigl)/4, \ldots$, and remembering that $\Ca_{ab}{}^{cd} = C_{ab}{}^{cd}$ is trace-free.

Multiplying by $V_a$   gives
$$   \eqalign{ & V_a  \Bigl(c_0  \CN^{a}{}_{c}+  c_2 \CNmm^{a}{}_{c}+c_3 \CNmmm^{a}{}_{c}+c_4 \CNmmmm^{a}{}_{c}+ \ldots  + c_{N-2} \Cb^{a}{}_{c} + {N-1\over 2} c_{N} \delta^{a}_{c}\Bigr) \cr &
= V_aC^{ia}{}_{pq} \Bigl(c_0 \CNm^{pq}{}_{ci}+ c_2 \CNmmm^{pq}{}_{ci}+c_3 \CNmmmm^{pq}{}_{ci}+c_4 \CNmmmmm^{pq}{}_{ci}+ \ldots  + c_{N-2} C^{pq}{}_{ci} \Bigr)+ {N-1\over 2} c_{N} V_c
  }$$

\

\

Multiplying by $V_a$   gives
$$ \eqalign{  0  & = V_a\Bigl(  c_0 \CN^{ab}{}_{cd}+ c_2 \CNmm^{ab}{}_{cd}+c_3 \CNmmm^{ab}{}_{cd}+c_4 \CNmmmm^{ab}{}_{cd}+   \cr & \qquad \qquad \ldots  + c_{N-2} \Cb^{ab}{}_{cd}  + c_{N-1} \Ca^{ab}{}_{cd}\Bigr) + c_{N} V_{[c}\delta^{b}_{d]}
\cr & 
 = V_a \Ca^{ab}{}_{ij}\Bigl(  c_0 \CNm^{ij}{}_{cd}+ c_2 \CNmmm^{ij}{}_{cd}+c_3 \CNmmmm^{ij}{}_{cd}+c_4 \CNmmmmm^{ij}{}_{cd}+ \cr & \qquad \qquad  \ldots  + c_{N-2} \Ca^{ij}{}_{cd} \Bigr) + c_{N-1} \Ca^{ij}{}_{cd} + c_{N} V_{[c}\delta^{b}_{d]}
 }$$

From which, by taking the trace, we obtain, 
$$ \eqalign{  0  & =  V_a C^{ab}{}_{ij}\Bigl(  c_0 \CNm^{ij}{}_{cb}+ c_2 \CNmmm^{ij}{}_{cb}+c_3 \CNmmmm^{ij}{}_{cb}+c_4 \CNmmmmm^{ij}{}_{cb}+ \cr & \qquad \qquad  \ldots  + c_{N-2} C^{ij}{}_{cb} \Bigr)  +{n-1\over 2} c_{N} V_d
 }$$
where $c_0 = 1, c_2 = -\Cb/2, c_3=  -\Cc/3, c_4= -\Bigl(\Cd-\Cb^2/2\Bigl)/4, \ldots$, and remembering that $ \Ca_{ab}{}^{cd} = C_{ab}{}^{cd}$, and so is trace-free.

\

\

It is important to realise that there is an important difference between the two sets of invariants $\{C^2, C^3, C^4, C^5 \}$ and $\{C^2, C^3, CC*,(C*)^3 \}$: although both sets are algebraically independent, the second set is also {\ algebraically complete} in the sense that remaining  Weyl scalar invariants can be constructed as polynomial functions of the four invariants, and none of these four invariants can be constructed as polynomial invariants of the  remainder; the first set does not have this property.

 \

 \
 
 \
 
 In the four-dimensional case, instead of using (10)  for $K^a$ (involving a product of six Weyl tensors divided by $c_6$ )   KNN  used (1) (involving a product of only two Weyl tensors divided by $c_22$), derived from the identity (4a). This certainly seems a considerable simplification; but on deeper analysis we note that 
 there are three alternative identities in [1] 
$$\eqalign{
\Cc^{cj}{}_{ci} & =  \delta^j_i C^3/4\cr
\Cd^{cj}{}_{ci} & =  \delta^j_i C^4/4\cr
\Ce^{cj}{}_{ci} & =  \delta^j_i C^5/4\cr} .\eqno(4b,c,d)$$
 Clearly we can solve for $K^a$ by this method providing there exists  at least one non-zero  Weyl scalar invariant from the set,  $\{C^2, C^3, C^4, C^5 \} $. 
 Taking these other possibilities into account, our version (10) for $K^a$ is not so unreasonably complicated, since the worst scenario case would require the use of the identity (4d) would give an expression for $K^a$ with five Weyl tensors,
 $$K^a = 4 \Cd^{ak}{}_{ij}\nabla^cC^{ij}{}_{kc}  /C^5 \eqno(11)
$$

 On the other hand, just as we found the Weyl scalar $C^6$ to be functionally dependent on lower Weyl scalar invariants, we would also expect there to be identities between the four index .....

 to be dependent on the lower powers of ${\bf C}$. These expections are confirmed by a little rearranging since  we can show
$$\Cf^a{}_b = 
$$ 

So therefore the expression 

From the four dimensional example we see that, with the use of  dimensionally dependent identities for each $n\ge 4$, we get lower orders of Weyl in the  forms for $K^a$ and a more useful statement of the restriction (10). So we now turn to the general case.

 at each order of Weyl invariants  there is at most one invariant which is independent of lower orders e.g. $C^a{}_i{}^b{}_j C^i{}_c{}^j{}_d C^c{}_a{}^d{}_b = 2 \Cc$ in four dimensions, but not in arbitrary dimensions;

\

\

\

The complex Weyl tensor  ${\cal C}^{ab}{}_{cd} = C^{ab}{}_{cd}+i C^{ab}{}_{cd}^*$  or equivalently  the Weyl spinor ${\bf \Psi} $ both lead to a tracefree {\it complex} $3\times 3$ matrix  ${\cal C}$ which incorporates the first Bianchi identity, and this representation is  simpler to manipulate than the $6\times 6$ matrix ${\bf C}$. We can now use the Cayley-Hamilton theorem for the complex trace-free matrix ${\cal C}$
$${\cal C}^3 -{1\over 2}<\!\!{\cal C}^2\!\!> {\cal C} -{1\over 3} <\!\!{\cal C}^3\!\!>{\bf I}=0
$$
to obtain relations between the dual invariants and the original $\Cp$ invariants. 

Beginning with 
$$<\!\!{\cal C}^4\!\!> -{1\over 2}<\!\!{\cal C}^2\!\!> <\!\!{\cal C}^2\!\!> =0
$$
we get 
$$2\Bigl(<\!\!\bigl(\Re({\cal C})\bigr)^4\!\!> + <\!\Im({\cal C})\bigr)^4\!\!>-6<\!\!\bigl(\Re ({\cal C}) \Im  ({\cal C})\bigr)^2\!\!>\Bigr)=\Bigl(<\!\! \bigl(\Re({\cal C})\bigr)^2 \!\!> -<\!\!\bigl(\Im({\cal C})\bigr)^2\!\!>\Bigr)^2 - 4 <\!\!(\Re({\cal C})\Im({\cal C})\!\!>^2
$$
and using the relationships between the real matrix ${\bf C}$ and the complex matrix ${\cal C}$
$${\bf C}=\pmatrix {{\Re( {\cal C})} & {\Im( {\cal C})} \cr- {\Im( {\cal C})} & {\Re( {\cal C})}}
$$
we get
$$<\!\!{\bf C}^2\!\!> = 2\Bigl(<\!\! \bigl(\Re({\cal C})\bigr)^2 \!\!> -<\!\!\bigl(\Im({\cal C})\bigr)^2\!\!>\Bigr)
$$
$$<\!\!{\bf C}^4\!\!> = 2\Bigl(<\!\!\bigl(\Re({\cal C})\bigr)^4\!\!> + <\!\Im({\cal C})\bigr)^4\!\!>-6<\!\!\bigl(\Re ({\cal C}) \Im  ({\cal C})\bigr)^2\!\!>\Bigr)
$$
It therefore follows directly that 
$$4 \Cd  = \Cb^2 - 16?(C^*C)^2\eqno(12)$$
and by similar calculations we also obtain 
$$\eqalign{
12 \Ce &  =5\bigl(\Cb\Cc- (C^*C)(C^*C^*C^*)\bigr)\cr
48\Cf & = 3\bigl(\Cb^3-3\Cb(C^*C)^2)\bigr)+8\bigl(\Cc^2-(C^*C^*C^*)^2\bigr)}\eqno(12)$$
From these expressions we are able to identify the relation between $\Cd,\Ce, \Cf$
$$  \eqalign {5\bigl(\Cb^2 & -4\Cd\bigr){\over 6!}\Bigl(24 \Cf- 18  \Cb \Cd
-8 \Cc^2 +3 \Cb^3\Bigr) \cr &
-4\Bigr(5\bigl(\Cb^2-4\Cd\bigr)\Cc^2  +{12\bigl(\Ce-\Cb\Cc\bigr)^2}\Bigr)=0}\eqno(13)$$
which  gives us an alternative expression to () for $\det ({\bf C})$
$$ \det ({\bf C}) = -20\Bigr(\Cc^2 +{12\bigl(\Ce-\Cb\Cc\bigr)^2\over 5\bigl(\Cb^2-4\Cd)}\Bigr)
\eqno(14)$$
and from this we conclude, providing $\Cb)^2\ne4\Cd$, 
$$\det ({\bf C})\ne 0 \quad \Leftrightarrow \quad \hbox{At least one of \Ce and \Cc} \  \ne 0 \ .
$$

\end